\documentstyle[12pt]{article}
\begin{document}
\noindent
\begin{center}
{\bf A CERTAIN CLASS OF LAPLACE TRANSFORMS WITH APPLICATIONS TO REACTION AND REACTION-DIFFUSION EQUATIONS}\\[0.5cm]
A.M. MATHAI \\
Department of Mathematics and Statistics, McGill University\\
Montreal, Canada H3A 2K6\\[0.5cm] 
R.K. SAXENA \\
Department of Mathematics and Statistics, Jai Narain Vyas University\\
Jodhpur-342005, India\\[0.5cm]
H.J. HAUBOLD\\
Office for Outer Space Affairs, United Nations\\
P.O.Box 500, A-1400 Vienna, Austria
\end{center}
\bigskip 
\noindent
{\bf Abstract.}
A  class of Laplace transforms is examined to show that particular cases of this class are associated with production-destruction and reaction-diffusion problems in physics, study of differences of independently distributed random variables and the concept of Laplacianness in statistics, $\alpha$-Laplace and Mittag-Leffler stochastic processes, the concepts of infinite divisibility and geometric infinite divisibility problems in probability theory and certain fractional integrals and fractional derivatives. A number of applications are pointed out with special reference to solutions of fractional reaction and reaction-diffusion equations and their generalizations. 

\section{Introduction}
We start with simple examples and then will proceed into generalizations. Consider a gamma density

  $$
f_1(x) = \frac{x^{\alpha -1} {\rm e}^{-\frac{x}{\beta}}}{\beta^\alpha \Gamma (\alpha )} ,~~ x > 0,~ \beta > 0, Re (\alpha ) > 0,
$$
where $Re (\cdot)$ denotes the real part of $(\cdot)$. Its Laplace transform, with parameter $p$, is given by the following:

\begin{equation}
L_{f_1} (p) = \int_0^\infty {\rm e}^{-px} f_1(x) {\rm d}x = (1+ \beta p)^{-\alpha}, 1+\beta p  > 0.
\end{equation}
For $\alpha = 1$ we have the Laplace transform of the exponential density. Consider the Laplace transform  of the Laplace density 

$$ f_2(x) = \frac{1}{2 \beta} {\rm e}^{-\frac{|x|}{\beta}}, -\infty < x < \infty,\;\; \beta > 0. $$
Its Laplace transform is given by the following

\begin{equation}
L_{f_2}(p) =  (1 - \beta^2 p^2)^{-1}.
\end{equation}
This Laplace transform $L_{f_2} (p)$ can also arise from a production-destruction or input-output model of  the following type: Consider the residual effect of an input-output type of situation $ u = x_1 - x_2$ where $  x_1$ and $ x_2$ are independently and identically distributed gamma type input and  gamma type output variables respectively. Then the Laplace transform of the density of  $u$, denoted by $L_u(p)$, is given by 

\begin{equation}
L_u (p) =  L_{x_1} (p) L_{x_2} (-p) = (1 + \beta p)^{- \alpha}  (1 - \beta p)^{- \alpha} =  (1 - \beta^2 p^2)^{- \alpha}.
\end{equation}
This residual variable $u$ has applications in many different areas. A few of the applications are mentioned in Mathai (1993a) and in Mathai, Provost, and Hayakawa (1995). The  Laplace transform  in (3) is also associated with the concept of Laplacianness in statistics, see Mathai (1993). This concept is connected to bilinear forms, quadratic form and the concept of chi-squaredness of quadratic forms, which is the basis for making inference in analysis of variance, analysis of covariance, regression and general model building areas.

If there are several independent input variables $x_1, \cdots, x_n$ such as the situation in  reaction or production problems,  and if there are several  independent output variables $x_{m+1}, \cdots, x_{m+n}$ and if they are all gamma type variables with different parameters then  the residual $u = x_1 + \cdots + x_m - x_{m+1} - \cdots - x_{m+n}$ has the Laplace transform

\begin{eqnarray}L_u (p) &=&   (1 + \beta_1 p)^{- \alpha_1} \cdots (1 + \beta_m p)^{- \alpha_m} \nonumber\\
& \times &(1 - \beta_{m+1} p)^{- \alpha_{m+1}} \cdots (1 - \beta_{m+n} p)^{- \alpha_{m+n}}.
\end{eqnarray}
The density of $u$, corresponding to the Laplace transform in (4), is also available in the literature, see for example Mathai and Provost (1992). Now consider the Mittag-Leffler function

$$
E_\alpha (x) =  \sum_{k=o}^\infty \frac{x^k}{\Gamma (1 + k \alpha)},~~ x > 0, Re (\alpha ) > 0,$$
and let
$$f_3(x)  =  E_\alpha (- x^\alpha),$$
then the Laplace transform
\begin{equation}
L_{f_3}(p) = \int_0^\infty {\rm e}^{-px} E_{\alpha} (-x^\alpha ) {\rm d}x = p^{\alpha -1} ( 1 +  p^\alpha )^{-1}.
\end{equation}

Haubold and Mathai (1995, 2000) considered a reaction problem where the number density of the reacting particles is a function of time $t$. For the $i^{th}$ particle let the number density be $N_i = N_i (t)$ with $ N_i (t = 0) = N_0$. If the production rate is proportional to the number density  then we have the simple differential equation, dropping $i$,

$$\frac{ {\rm d}}{ {\rm d}t} N (t) = a~~  N (t),~~~$$
where $a$ is a constant. If some particles produced are also destroyed or consumed and if the destruction rate is given by

$$\frac{ {\rm d}}{ {\rm d}t} N (t) = -b ~~ N (t),~~~$$
then the residual effect is given by 
$$\quad\quad~~~\frac{ {\rm d}}{ {\rm d}t} N (t) = -c~~  N (t), \; c = b - a, $$
and then the solution is  
$$ N (t) = N_0 {\rm e}^{-ct}.$$
But if a fractional integral, instead of the full integral, is used then the reaction equation is given by

\begin{equation}
N (t) - N_0 = -c^\nu {_0D_t}^{-\nu} N(t), \nu > 0,
\end{equation}
where $c$ is replaced by $c^\nu$  and its solution has the form

$$ N (t) = N_0 \sum_{k=0}^\infty \frac{(-1)^k(ct)^{k\nu }}{\Gamma (1 + k \nu )}=  N_0 E_\nu (-c^\nu t^\nu ),$$
where the fractional integral  is defined as follows

\begin{equation}
{_aD_t}^{-\nu} f(t) = \frac{1}{\Gamma ( \nu )} \int_a^t (t-u)^{\nu-1} f(u){\rm d}u,\;\;\nu > 0,
\end{equation} 
with ${_aD_t}^{0}f(t) = f(t)$. The Laplace transform of $N(t)$, coming from (6), is then,

\begin{equation}
L_{ N (t)} (p) = \frac{N_0}{p[1+(\frac{c}{p})^\nu ]},
\end{equation}
which is a constant multiple of (5) for $c=1, \nu = \alpha$.

\section{Generalizations}
If the basic fractional production-destruction equation in (6) is modified to the form 

\begin{equation}
N(t) - N_0 t^{\mu - 1} = -c^\nu {_0D_t}^{-\nu}N(t), \mu > 0, \nu > 0,
\end{equation}
then the Laplace transform is given by 
\begin{eqnarray}
L_{ N (t)} (p) &=& \frac{N_0 \Gamma (\mu)}{p^\nu[1+(\frac{c}{p})^\nu ]} = \frac{N_0 \Gamma (\mu)}{p^{\mu - \nu} [{p}^\nu+ {c}^\nu]}\\
& = &\frac{N_0 \Gamma (\nu)}{[{p}^\nu+ {c}^\nu]}\;\;  \mbox{for}\;\;  \mu = \nu. 
\end{eqnarray}
This Laplace transform in (11) for $ c = 1$ is associated with infinitely divisible and geometrically infinitely divisible distributions, $\alpha$-Laplace and Linnik distributions in probability theory and $\alpha$-Laplace and Mittag-Leffler stochastic processes, see the works of R.N. Pillai and his associates, a summary of which is available from Jose and Seetha Lekshmi (2004). A comprehensive presentation of this field of processes is given in Uchaikin and Zolotarev (1999). The class of Laplace transforms relevant in geometrically infinitely divisible and $\alpha$-Laplace distributions is of the form 

\begin{equation}
L(p) = \frac{ 1}{1+ \eta (p)},
\end{equation}
where $\eta (p)$ satisfies the condition $\eta (bp) = b^\alpha \eta (p)$ where $\eta (p)$ is a periodic function for fixed $\alpha $. Observe that (11) satisfies the above condition with $\eta (p) = p^\nu$ and $\alpha = \nu$. Thus the basic reaction equation in (9) is also associated with certain stochastic processes and infinitely divisible distributions.

Another observation that can be made is that the inverse of (10) or the number density $N(t)$, which gives rise to (10), is of the following general form

\begin{equation}
N (t) = N_0 \Gamma (\mu) t^{\mu -1} E_{\nu, \mu} (- c^\nu t^\nu), 
\end{equation}
where $ E_{\nu, \mu}(\cdot )$ is a generalized Mittag-Leffler function,

\begin{equation}
E_{\nu, \mu} (-c^\nu t^\nu ) =  \sum_{k=0}^\infty \frac{(-1)^k(t^\nu c^\nu)^{k }}{\Gamma (\mu  + k \nu )},~ \mu > 0,~ \nu > 0.
\end{equation}
Note that (14) is also connected to the F-function introduced by Hartley and Lorenzo (1998) and Lorenzo and Hartley (1999). Their F-function is the following,

\begin{equation}
F_q (-a, t) =  t^{q-1} \sum_{n=0}^\infty \frac{(-a)^n~  t^{n q}}{\Gamma (q + n q )} =  t^{q-1} E_{q,q} (-a t^q ), Re(q) > 0.
\end{equation}
When the variable is restricted to the interval $ t \geq \delta$ for some $\delta ,$ then (15) reduces to the R-function of Lorenzo and Hartley (1999), which is defined as

\begin{eqnarray}
R_{\nu , \mu} (a, \delta , t) &=& \sum_{n=0}^\infty \frac{a^n (t - \delta )^{(n  + 1)\nu - \mu -1}}{\Gamma [(n+1) \nu -\mu]},~~ t > \delta > 0\nonumber\\
&=& (t - \delta )^{\nu -  \mu - 1} E_{\nu , \nu -\mu }[ a(t - \delta)^\nu ],~~ t > \delta > 0.
\end{eqnarray}
If we consider another modification of the basic reaction equation in (9) to the form

\begin{equation}
N(t) - N_0 ~ t^{\mu -1} E_{\nu, \mu}^\gamma (-c^\nu t^\nu) = -c^\nu {_0D_t}^{-\nu} N(t), 
\end{equation}
where $E_{\nu, \mu}^\gamma (\cdot ) $ is a further generalized form of the Mittag-Leffler function given by 

\begin{equation}
E_{\nu, \mu}^\gamma (x) = \sum_{k=0}^\infty \frac{(\gamma)_k}{k!} \frac{x^k}
{\Gamma (\mu +k\nu)} ~ Re (\mu) > 0, ~ Re (\nu ) > 0, 
\end{equation}
$$(\gamma)_k = \gamma (\gamma +1) \cdots  (\gamma +k-1), ~ \gamma \ne 0, (\gamma)_0 = 1,$$
then the Laplace transform is the following

\begin{equation}
 L_{N(t)} (p) = \frac{ N_0}{p^{\mu - \nu (\gamma +1)} [c^\nu + p^\nu ]^{\gamma +1}} 
.\end{equation}
The inverse of (19) or the number density in this case is given by

\begin{equation}
N(t) =  N_0 ~ t^{\mu -1} E_{\nu ,  \mu }^{\gamma +1} (-c^\nu t^\nu). \end{equation}
The most general form of Laplace transform associated with Mittag-Leffler functions is that in (19), which is a special case of the general class of Laplace transforms associated with $\alpha$-Laplace stochastic processes and geometrically infinitely divisible statistical distributions.

If the Mittag-Leffler function is further generalized then the Laplace
transforms will enter the family of Fox's H-function. For example,
let us consider a simple generalization of the Mittag-Leffler function
to the following form
        $$ g_1(t) = \sum_{k=0}^\infty \frac{ (\gamma_1)_k}{(\beta_{1})_k k!}~~ \frac{t^k}{\Gamma (\beta + k\alpha)}.$$
        Then the   Laplace transform of $t^{\beta -1}g_1(t^\alpha)$,   denoted by $L_1(p)$, is given by the following

       \begin{equation}
       L_1(p) = \sum_{k=0}^\infty \frac{ (\gamma_1)_k}{(\beta_1)_k k!}~~ p^{-(\alpha k + \beta ) }
        =p^{-\beta} {_1F_1 }(\gamma_1; \beta_1 ; p^{- \alpha}),
        \end{equation}
        where $ {_1F_1 }$ is a confluent hypergeometric function. If the factor
        $t^{\beta-1}$ is not incorporated then the Laplace transform will be a special case
        of an H-function. For a description and properties of  Fox's H-function see Mathai and Saxena
        (1978), Mathai (1993a), and Kilbas and Saigo (2004).

        A direct generalization of a Mittag-Leffler  function  and the most generalized form in this category is
        Wright's function given by

         \begin{equation}
         {_p\psi_q} \biggl[_{(b_1, B_1), \cdots,  (b_q, B_q)}^{(a_1, A_1), \cdots,  (a_p, A_p)} ; z\biggr] =
           \sum_{k=0}^\infty \frac{\left\{ \prod_{j=1}^p \Gamma (a_j +A_jk)\right\} }
           {\left\{ \prod_{j=1}^q \Gamma (b_j +B_j k)\right\} } \frac{z^k}{k!},
           \end{equation}
           with  $1+ \sum_{j=1}^q B_j - \sum_{j=1}^p A_j \ge 0$.\\

\noindent
           Wright's function is again a particular case of Fox's H-function. Note that the
           Mittag-Leffler function

           \begin{equation}
        E_{\alpha , \beta }(z) =  {_1\psi_1} \biggl[_{(\beta , \alpha)}^{(1, 1)} ; z\biggr] =
        H_{1 ,2}^{1 ,1} \biggl[ -z\big\vert_{(0 ,1), (1- \beta , \alpha)}^{(0,1)}\biggr],
        \end{equation}
        where $ H(\cdot)$ denotes an H-function, and the H-function is defined in terms of
        the following Mellin-Barnes type representation

         \begin{equation}
         H_{p , q}^{m , n} \biggl[ z \big\vert_{(b_1, B_1), \cdots,  (b_q, B_q)}^{(a_1, A_1), \cdots,  (a_p, A_p)} \biggr]
         = \frac {1}{2 \pi i} \int_L g(s) z^{-s} {\rm d}s, i = \sqrt{-1},
         \end{equation}
         where

\begin{equation}
         g(s) = \frac{\left\{ \prod_{j=1}^m \Gamma (b_j +B_j s)\right\}  {\left\{ \prod_{j=1}^n \Gamma (1-a_j - A_j s)\right\} }}
         {\left\{ \prod_{j=m+1}^q \Gamma (1-b_j - B_j s)\right\}  {\left\{ \prod_{j=n+1}^p \Gamma (a_j +A_j s)\right\}}},
          \end{equation}
          with $A_j, j=1, \cdots , p$ and     $B_j, j=1, \cdots , q$ being positive real numbers and  $a_j$'s and $b_j$'s are complex quantities. Existence conditions, the types of contour $L$ and properties are available from Mathai and Saxena (1978) and Kilbas and Saigo (2004).
        
        Observe that for $\gamma = 1$ in (18) we have $E_{\nu ,
          \mu}^1 = E_{\nu,\mu }$. Hence for  \mbox{$\gamma = 1,2,3\cdots $,}
        for positive integer values, one can write $E_{\nu ,
          \mu}^\gamma$ as a linear function of $E.,.(\cdot )$. Thus,
        certain linear functions of Laplace transforms of the type in
        (10) can be summed up to a Laplace transform of the type in
        (19). Some examples of this type may be seen from a series
        of papers by Haubold, Mathai, and Saxena, for example, see
        Saxena, Mathai, and Haubold (2004).  One such example can be
        stated as follows: If $c > 0, \mu > 0, \nu > 0$ then the
        solution of the fractional integral equation for the number
        density $N(t)$,

\begin{equation}
N(t) - N_0 ~ t^{\mu -1} E_{\nu, \mu }(- c^\nu t^\nu ) = -c^\nu {_0D_t}^{-\nu} N(t),
          \end{equation}
          is given by

          \begin{equation}
          N(t) = \frac { N_0~ t^{\mu -1}}{\nu}\biggl[ E_{\nu, \mu -1}(-c^\nu t^\nu ) +(1 - \mu+\nu)
          E_{\nu, \mu} (- c^\nu t^\nu)\biggr].
          \end{equation}

\section{Reaction Equation}
If the Laplace transform has the structure of a product of
          various factors then such a case can also be handled without
          much difficulty.  As an example, consider the production-destruction
          fractional integral model for the number density $N(t)$ of
          the following type

          \begin{equation}
          N(t) -  N_0 ~t^{\mu -1}  E_{\nu, \mu }(- d^\nu t^\nu )  = - c^\nu {_0D_t}^{-\nu} N(t).
 \end{equation}
         Then it is not difficult to see that the Laplace transform
          of $N(t)$, with parameter $p$, denoted by $L_2(p)$ is the
          following

   \begin{equation}
  L_2(p)  = \frac{N_0} {{p^\mu } [ 1+ (\frac{c}{p})^\nu ]
                 [ 1+ (\frac{d}{p})^\nu ]}.
             \end{equation}
             When $c = d$ then this case reduces to the Laplace
             transform of a generalized Mittag-Leffler function,  which can also be
             written as a linear combination of simple Mittag-Leffler
             functions as explained in (26) and (27). When $c
             \neq d$ then we may consider the following identity

\begin{equation}
 \frac{1}{(p^\nu + c^\nu) (p^\nu + d^\nu) } = \frac{1}{(c^\nu -d^\nu)}\biggl[ \frac{1}{p^\nu + d^\nu} - \frac{1}{p^\nu + c^\nu }  \biggr], c \neq d.            \end{equation}
 Hence,
 
\begin{equation}
      L_2(p) =  \frac{1}{(c^\nu - d^\nu)} \left\{ \frac{1} {p^{\mu -\nu} }~  \frac{1}{1+ (\frac{d}{p})^\nu}
     -  \frac{1}{p^{\mu -\nu} } ~ \frac{1}{1+ (\frac{c}{p})^\nu }\right\}.
      \end{equation}
Then taking the inverse Laplace transform, we have,
\begin{equation}
N(t) = N_0 ~ \frac{t^{\mu - \nu -1}}{(c^\nu - d^\nu ) }\biggl[E_{\nu , \mu -\nu} (- {\rm d}^\nu t^\nu) - E_{\nu , \mu -\nu} (-c^\nu t^\nu )\biggr] .
\end{equation}

Now, let us look into another situation where the Laplace transforms are of the following types

\begin{equation}
L_3(p) = \frac{p^{\alpha -1}}{p^\alpha + a p^\beta + b}\;\mbox{and}\; L_4(p) =
\frac{p^{\beta -1}}{p^\alpha + a p^\beta + b},
\end{equation}
where $a, b , \alpha , \beta $ are constants. Such Laplace transforms can be
 handled by using the following procedure. For example, consider

\begin{eqnarray}
L_3(p) &=&   \frac{p^{\alpha -1}}{p^\alpha + a p^\beta + b} = \frac{p^{\alpha - \beta -1}}{(p^{\alpha - \beta}+ b p^{-\beta}) [ 1+ \frac{a}{(p^{\alpha - \beta}+ b p^{-\beta})}]}\nonumber\\
&=& \sum_{r=0}^\infty (-a)^r \frac{p^{-(\alpha - \beta ) r-1}}{(1+ \frac{b}{p^\alpha})^{r+1}}.
\end{eqnarray}
Comparing with the Laplace transform of a generalized Mittag-Leffler function we have the inverse
given by the following

\begin{equation}
L^{-1} \biggl[ \frac{p^{\alpha -1}}{p^\alpha + a p^\beta + b} \biggr]
= \sum_{r=0}^\infty (-a)^r t^{(\alpha - \beta ) r } E_{\alpha, (\alpha - \beta)
r+1}^{r+1} (-bt^\alpha).
\end{equation}

\section{Reaction-Diffusion Equation}
A problem recently considered by Saxena, Mathai, and Haubold (2005) is a reaction-diffusion system occurring in various areas of pattern formation in biology, chemistry, physics, see for example Henry and Wearne (2000, 2002), Henry, Langlands, and Wearne (2005), and Manne, Hurd, and Kenkre (2000). The basic reaction-diffusion equation has the form

\begin{equation}
\frac{\partial N}{\partial t} = \mu \frac{\partial^2  N}{\partial x^2} + \lambda\  f(N), N = N(x,t),
\end{equation}
where $\mu$ is the diffusion constant, $\lambda $ is a constant and $f(N)$ is a
nonlinear function of $N$. When $f(N) = \delta N (1-N)$, where $\delta$ is a
constant, then (36) reduces to the real Ginsburg-Landau equation. An equation in this category examined by
Manne, Hurd, and Kenkre (2000) is the following

\begin{equation}
\frac{\partial^2 N}{\partial t^2} + a \frac{\partial N}{\partial t}  = \nu^2 \frac{\partial^2  N}{\partial x^2}
\xi^2 N (x,t),
\end{equation}
where $\xi$ indicates the strength of nonlinearity in the system. For
solving (37) one can adopt the procedure of taking the Laplace
transform with respect to $t$ and Fourier transform with respect to
$x,$ simplifying and then inverting the Laplace-Fourier transform to
obtain $N(x,t)$. Details of the procedure may be seen from Saxena, Mathai, 
and Haubold (2004, 2005). In this procedure, the crucial step is the inversion of the
Laplace transform of the type $L_4 (p) $ in (33). From (35) one can
see that the inverse can be written as a series of Mittag-Leffler
functions

\begin{equation}
L^{-1} \biggl[ \frac{p^{\beta -1}}{p^\alpha + a p^\beta + b} \biggr]
= \sum_{r=o}^\infty (-a)^r t^{(\alpha - \beta ) (r+1) } E_{\alpha, (\alpha - \beta)
(r+1)+1}^{r+1} (-bt^\alpha).
\end{equation}

\section{Conclusions}
Linear and nonlinear reaction, diffusion, and reaction-diffusion equations, respectively, are used to model spatio-temporal processes in physical systems, for example astrophysical fusion plasmas (Wilhelmsson and Lazzaro, 2001; Kulsrud, 2005) . Such equations give rise to dissipative structures and self-organization phenomena (Nicolis and Prigogine, 1977; Haken 2004). Attempts have been made to employ them to discover a time-dependent mechanism in solar nucleosynthesis (Haubold and Mathai, 1995, 2000). 

The linear reaction (relaxation) equation can be used to explore fundamental principles of Boltzmann-Gibbs statistical mechanics and its nonlinear generalization leads to new insights into nonextensive statistical mechanics (Tsallis, 2004). This approach has also been extended to discover the connection between Tsallis' nonextensive maximum entropy formalism and the nonlinear reaction-diffusion equation (Tsallis and Bukman, 1996).

In this paper we used Laplace transform techniques to derive closed-form solutions of reaction equations starting from its simplest form and its fractional version in eq. (6) through generalizations of the equation given in (9), (17), (26), and (28). Along this way of generalizing the linear reaction equation we provide their respective closed-form representations in terms of Mittag-Leffler functions. Gradually the link between the Mittag-Leffler function and Wright's function and, subsequently, Fox's function has been shown. This transition can be characterized as a transition from the standard exponential behavior ("normal" relaxation) to a  power-law behavior ("anomalous" relaxation) of the solutions of the respective reaction equations. For the reaction-diffusion equations (36) and (37) we indicate that the same Laplace transform techniques, as successfully employed for deriving closed-form solutions of reaction equations, can also be used for finding solutions of reaction-diffusion equations.\par  
\bigskip
\noindent
{\bf Acknowledgment}\\
A.M.M. would like to thank the Natural Sciences and the Engineering
Research Council of Canada for the grant in carrying out this research.\par

\section*{References}\par
\noindent
Haken, H.: 2004, \emph {Synergetics: Introduction and Advanced Topics}, 

Springer-Verlag, Berlin-Heidelberg.\\
Hartley, T.T. and Lorenzo, C.F.: 1998, A solution to the fundamental 

linear fractional order differential equation, \emph{NASA/TP-1998-208693}.\\
Haubold, H.J. and Mathai, A.M.: 1995, A heuristic remark on the periodic 

variation in the number of solar neutrinos detected on Earth,

\emph {Astrophysics and Space Science} {\bf 228}, 113-134.\\ 
Haubold, H.J. and Mathai, A.M.: 2000, The fractional kinetic equation 

and thermonuclear functions, \emph {Astrophysics and Space Science} 

{\bf 327}, 53-63.\\
Henry, B.I. and Wearne, S.L.: 2000, Fractional reaction-diffusion,

\emph{Physica A} {\bf 276}, 448-455.\\
Henry, B.I. and Wearne, S.L.: 2002, Existence of Turing instabilities 

in a two-species fractional reaction-diffusion system,
 
\emph {SIAM Journal of Applied Mathematics} {\bf 62}, 870-887.\\
Henry, B.I., Langlands, T.A.M., and Wearne, S.L.: 2005, Turing pattern 

formation in fractional activator-inhibitor systems, \emph {Physical Review E} 

{\bf 72}, 026101.\\
Jose, K.K. and Seetha Lekshmi, V. 2004, \emph{Geometric Stable Distributions:}

\emph{Theory and Applications}, Science Educational Trust Pala.\\
Kilbas, A.A. and Saigo, M.: 2004, \emph{H-Transforms: Theory and Applications}, 

Chapmanand Hall/CRC, Boca Raton-London-New York-Washington, D.C.\\
Kulsrud, R.M.: 2005, \emph {Plasma Physics for Astrophysics}, Princeton 

University Press, Princeton and Oxford.\\
Lorenzo, C.E. and Hartley, T.T.: 1999, Generalized functions for the 

fractional calculus, \emph {NASA/TP-1999-209424/REV1}.\\
Manne, K.K., Hurd, A.J., and Kenkre, V.M.: 2000, Nonlinear waves in 

reaction-diffusion system: The effect of transport memory,  

\emph{Physical Review E} {\bf 61,} 4177-4184.\\
Mathai, A.M.: 1993, On non-central generalized Laplacianness of quadratic 
 
forms in random variables, \emph{ Journal of Multivariate Analysis}

{\bf 45}, 239-246.\\
Mathai, A.M.: 1993a, \emph{A Handbook of Generalized Special Functions}

\emph{for Statistics and Physical Sciences}, Oxford University Press, Oxford.\\
Mathai, A.M. and Provost, S.B.: 1992, \emph{Quadratic Forms in Random}

\emph{Variables:} \emph{Theory and Applications}, Marcel Dekker, New York.\\
Mathai, A.M.,  Provost, S.B., and Hayakawa, T.: 1995, \emph{Bilinear Forms} 

\emph{and Zonal Polynomials. Lecture Notes in Statistics}, Springer-Verlag,

New York.\\
Mathai, A.M. and Saxena, R.K.: 1978, \emph{The H-function with Applications} 

\emph{in Statistics and Other Disciplines}, Wiley Eastern, New Delhi and 

Wiley Halsted, New York.\\
Nicolis, G. and Prigogine, I.: 1977, \emph {Self-Organization in Nonequilibrium} 

\emph {Systems: From Dissipative Structures to Order Through Fluctuations}, 

John Wiley and Sons, New York.\\  
Saxena, R.K., Mathai, A.M., and Haubold, H.J.: 2004, On generalized

fractional kinetic equations, \emph{Physica A} {\bf 344}, 657-664.\\
Saxena, R.K., Mathai, A.M., and Haubold, H.J.: 2005, Reaction-diffusion

systems and nonlinear waves, this volume.\\
Tsallis, C.: 2004, What should a statistical mechanics satisfy to reflect

nature?, \emph {Physica D} {\bf 193}, 3-34.\\
Tsallis, C. and Bukman, D.J.: 1996, Anomalous diffusion in the presence 

of external forces: Exact time-dependent solutions and their 

thermostatistical basis, \emph {Physical Review E} {\bf 54}, R2197-R2200.\\
Uchaikin, V.V. and Zolotarev, V.M.: 1999, \emph{Chance and Stability:}

\emph{Stable Distributions and Their Applications}, VSP, Utrecht, 

The Netherlands.\\
Wilhelmsson, H. and Lazzaro, E.: 2001, \emph {Reaction-Diffusion Problems} 

\emph {in the Physics of Hot Plasmas}, Institute of Physics Publishing, 

Bristol and Philadelphia.\\
          
\end{document}